\documentstyle[amscd,amssymb,graphicx,12pt]{amsart}
\input xypic \xyoption{curve}
\input epsf

\def\hcorrection#1{\advance\hoffset by #1 }
\def\vcorrection#1{\advance\voffset by #1 }

\vcorrection{-.7in}
\hcorrection{-.2in} 

\newcommand{\C}[1]{{\cal#1}} 
\newcommand{\D}[1]{{\Bbb#1}} 

\theoremstyle{plain}
\newtheorem{th}{Theorem}[section]
\newtheorem{cor}{Corollary}[section]
\newtheorem{lem}{Lemma}[section]
\newtheorem{prop}{Proposition}[section]
\newtheorem{cl}{Claim}[section]
\newtheorem{conj}{The Coefficient Theorem}[section]

\theoremstyle{definition}
\newtheorem{defin}{Definition}[section]

\theoremstyle{definition}

\theoremstyle{remark}
\newtheorem{rem}{Remark}[section]

\numberwithin{equation}{section}

\def\la{\rightharpoonup}
\def\ra{\rightharpoondown}
\def\tg{\tilde{\Gamma}}
\def \n{\noindent }
\def \bs{\bigskip}

\begin{document}

\pagestyle{plain}
\addtolength{\footskip}{.3in}

\title[On coefficients of star-products]
{A canonical semi-classical star-product}
\author{Lucian M. Ionescu and Papa A. Sissokho}
\address{Department of Mathematics, Illinois State University, IL 61790-4520}
\email{lmiones@@ilstu.edu,psissok@@ilstu.edu} 

\begin{abstract}
We study the Maurer-Cartan equation of the pre-Lie algebra of graphs
controling the deformation theory of associative algebras
and prove that there is a canonical solution within the class of graphs without circuits,
without assuming the Jacobi identity.
The proof is based on the unique factorization property of graph insertions.
\end{abstract}

\maketitle
\tableofcontents


\section{Introduction}\label{S:intro}
In \cite{comb} it was claimed that the initial value deformation
problem in the pre-Lie algebra of graphs has a canonical solution when restricted to 
graphs without circuits.
The existence relied on Kontsevich solution, 
i.e. a star-product corresponding to a general Poisson structure,
which conjecturally yields a star-product when restricted to graphs without circuits (see also \cite{FI}).

In the case of linear Poisson structures 
star products have been given by S. Gutt \cite{Gutt,Gutt2} and 
studied in the light of Kontsevich approach by Polyak \cite{Po}.

In this article we investigate solutions of the Maurer-Cartan equation 
in a differential graded Lie algebra of graphs
from the combinatorial point of view
and study the pre-Lie algebra of the corresponding graphs.

The main result is an explicit solution of the Maurer-Cartan equation
in the differential graded Lie algebra of graphs which controlls
the deformation theory of associative algebras (Theorem \ref{T:main}).
The proof relies on the unique factorization of graph insertion at the level
of a boundary point (Corollary \ref{C:UF}).
Together with a result regarding the multiplicity coefficient for the 
the above mentioned graph insertion (The Coefficient Theorem \ref{conj}),
we prove that the ``graph exponential'' $\sum_\Gamma \Gamma/|Aut(\Gamma)|$
is a solution.
As a corollary, the cohomological obstructions vanish.
In particular the cohomology class of the Jacobiator is zero.

The article is organized as follows.
In Section \ref{S:pl} we recall the class of graphs \cite{Kon1} 
together with the pre-Lie composition from \cite{comb} (see also \cite{Po,Kath}).
The core of the article is Section \ref{S:solution}
which claims the ``obvious solution'' and introduces the main properties
of graph insertion used in the proof.
In Section \ref{S:conclusions} we discuss some related questions.

\section{The pre-Lie algebra of graphs}\label{S:pl}
The combinatorial problem regarding the coefficients of a star-product
is captured by the ``graphical calculus'' we will call {\em Kontsevich rule},
a sort of a ``dual Feynman rule'':
$$B:k\C{G}\to D, \quad B(\Gamma)=\C{U}(exp(\alpha)).$$
Here $\C{G}$ is a class of graphs,
$(D,\circ,m)$ is some {\em pointed pre-Lie algebra} \cite{HDGLA} 
with a distinguished element $m$ such that $m\circ m=0$
and $\alpha$ is a Poisson structure (say on $R^n$):
$$\alpha=\sum_i \alpha^{ij}\partial_i\otimes \partial_j.$$
It is an antisymmetric 2-tensor satisfying Jacobi identity:
$$\alpha^{ij}=-\alpha^{ji}, \quad \sum_{circular}\{\{f,g\},h\}=0.$$
To a particular type of Poisson structure (e.g. constant/linear coefficients)
corresponds a specific class of directed labeled graphs:
those graphs $\Gamma$ which are not in the kernel of the Kontsevich rule
($k\C{G}/Ker K$).
Once the ``Problem'' is pull back to graphs,
it amounts to solving the equation $Z\circ Z=0$:
$$Z=\sum_n Z_n\ h^n, \quad Z_n=\sum_{\Gamma\in \C{G}_{n,2}} W_\Gamma\ \Gamma,$$
in a pre-Lie algebra with composition $\circ$ \cite{Po},
defined independently in \cite{comb}.

%
%
\subsection{Lie admissible graphs}
Let $\tilde{\C{G}}_{n,m}$ be the set of {\em orientation classes} 
of {\em Lie admissible edge labeled graphs} of \cite{Po}, p.3,
corresponding to {\bf linear Poisson structures} (see also \cite{comb}).
An element $\tilde{\Gamma}\in \tilde{\C{G}}_{n,m}$ is a directed graph with $n$ internal vertices, 
$m$ labeled boundary vertices $1,2...,m$ such that 
each internal vertex is trivalent with exactly two descendants.
The corresponding arrows will be labeled left/right,
defining the {\em orientation class} of the graph $\Gamma$ 
up to a ``negation'' of the edge labeling in any two internal vertices \cite{Po}.
The corresponding (graded) is denoted by $\C{G}=\cup \C{G}_m$,
where $\C{G}_m=\cup_{n\in \D{N}} \C{G}_{n,m}$.

%
%
\subsection{Graphical representation and notation}\label{S:graphrep}
The order of the boundary vertices is ``fixed'' once and for all
and will be represented graphically by placing the boundary vertices 
on a oriented line.

The left/right labels on the outgoing edges at each internal vertex
are implicit in a {\bf graphical representation} of a graph $\tilde{\Gamma}$,
which is an {\em embedding} $\sigma:\tilde{\Gamma}\to H$ 
of the graph as a ``discrete manifold''
into the upper half-plane $H$ ($\partial H$ is the above oriented line) with some
metric such that each oriented pair of points determines a unique connecting geodesic.

The embedding is 1:1 at the level of internal vertices and maps the two outgoing ``tangent vectors''
at an internal vertex to a base of the tangent space at the corresponding point:
$$\forall v\in\tilde{\Gamma}^{(0)}\quad T\tau_v:T_v(\Gamma)\to T_{\tau(v)}(H)\quad isomorphism.$$
The left/right labeling of the arrows of $\Gamma^{(1)}$
is induce by the counterclockwise orientation of the plane ($H$)
such that $T\tau$ is an orientation preserving embedding.
In particular, the {\em outgoing angle} at an internal vertex embedded in $H$ is not $\pi$
(e.g. the embedding must ``brake the symmetry'' of graph $c_2$).

In what follows we will use ($\la$) to denote an arrow with an $L$ label (first vector)
and ($\ra$) for an arrow with and $R$ label (second vector). 
Moreover, we will use full arrows ($\rightarrow$) to denote the representative $[\Gamma]$ of 
all $graphs$ $\Gamma$ differing only by the labeling of their edges.
The corresponding graded set is denoted by $[\C{G}]=\cup [\C{G}]_{n,m}$.

The graphs from $\C{G}_{n,2}$ with $n=0,1,2$ internal vertices are 
{\em prime Bernoulli graphs} $b_0, b_1, b_2^{L/R}$
or the products of Bernoulli graphs ($b_1^2$; to be defined shortly).

\vspace{.2in}
\hspace{0.75in} 
\epsfbox{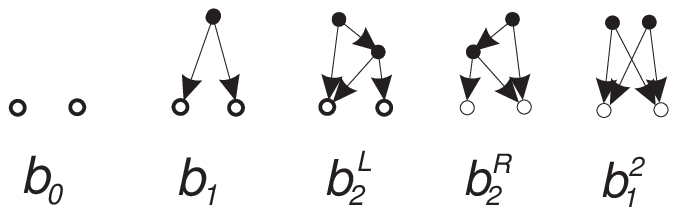} 

The prime graphs from $\C{G}_{n,3}$ are represented bellow.

\hspace{.2in} \epsfbox{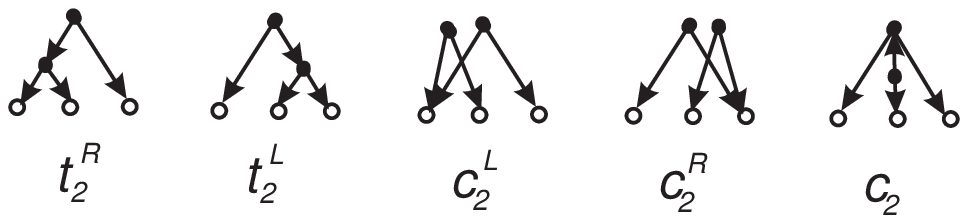}

%
%
\subsection{Antisymmetry relation}
The Kontsevich rule has an obvious kernel since the Poisson tensor $\alpha$
is antisymmetric and satisfies the Jacobi identity.
We will consider the two corresponding relations on graphs at distinct levels.

Denote by $H=k\C{G}_{n,m}/\sim$ the quotient modulo the equivalence relation
generated by $\Gamma'\sim -\Gamma$,
where $\Gamma'$ is obtained from $\Gamma$ by ``negation'' (switching) of the left/right labeling at 
one internal vertex only.
This completes the process of taking orientation classes of edge-labeled graphs
and will be considered independently of the Jacobi identity (compare \cite{Po}, p.5).
The same notation will be used for the induced quotient (linear) map $K:H\to D$.

Note at this point that $K$ cannot be defined on graphs after forgetting the edge-labels
(linear map $[\ ]$):
$$\diagram
\tilde{\C{G}} \drto_{ng} \ar@{->>}[r]^{ng^2} & \C{G} \rto^{K}\ar@{->>}[d]_{or} & D\\
& H \urto^{K} \ar@<2pt>[r]^{[\ ]} & \ar@<2pt>[l]^{\tau} k[\C{G}]\ar@{.>}[u]_{No}.
\enddiagram$$
Nevertheless the compositions of unlabeled graphs of \cite{comb} 
will be used in computing the compositions in $H$.
In order to do this a section $\tau$ may be defined by choosing an embedding
of each graph $\Gamma\in [\C{G}]$ (see \S\ref{S:graphrep}).
Additional procedures compatible with the additional structure will be considered later on.

Note that $[\ ]$ is a 2:1 covering map and $\C{G}\cong [\C{G}]\times \D{Z}_2$.
Moreover, with $k$ denoting the number of orientation inversions of a (edge) labeled graph, 
the following equality holds in $H$:
$$\forall \tilde{\Gamma}\in \tilde{\C{G}} \quad \tilde{\Gamma}=(-1)^k\tau([\Gamma]) \quad 
(\Gamma=ng(\tilde{\Gamma})).$$

%
%
\subsection{Product of graphs}
The {\em product of graphs} $\C{G}_{n,m}\times \C{G}_{n',m}\to \C{G}_{n+n',m}$
(L-graph multiplication \cite{Kath}, p.23; \cite{Po}, p.3, \cite{comb}, p.5)
is defined by identifying their corresponding boundary points.
A graph is {\em prime} if by ``cutting its boundary'' 
it yields a ``graph'' with only one component.
For example, $b_1^2=(b_1)^2$, is not a prime graph.

Note that the product is compatible with the equivalence relation on edge-labeled graphs,
inducing a product on $H$.

The subspace generated by prime graphs $\C{G}_\C{P}$ is denoted by $g$.
For this purpose the unit $b_0$ is considered prime.
\begin{prop}\label{P:tau1}
$H=k[\C{G}_\C{P}]=S^\bullet(g)$ is the polynomial algebra generated by the prime graphs.
Any section $\tau$ defined on prime graphs extends uniquely to $H$
as an algebra morphism.
\end{prop}

%
%
\subsection{Composition of graphs}
The graph composition of \cite{Kath} was introduced in \cite{comb}, p.8
at the level of unlabeled graphs,
as the pullback of the Gerstenhaber composition through Kontsevich rule (see also \cite{Po}).
It acquires Leibniz rule in this process,
since under Kontsevich representation the arrows carry differential operators,
while boundary vertices are ``colored'' by functions.
For example:
$$b_1^2\circ b_0=b_1^2\circ_1 b_0-b_1^2\circ_2 b_0=
\bullet b_1^2-b_1^2\bullet+2(c_2^R-c_2^L),$$
where with $\bullet$ the only graph in $\C{G}_{0,1}$
and $\bullet \Gamma$ denotes the ``concatenation'' of 
the corresponding graphs.

Note that graph composition is compatible with the grading by the
number of boundary vertices (see \cite{comb}, Appendix p.22):
$$\Gamma\in\C{G}_{n,m}, \quad deg_b(\Gamma)=m-1,$$
$$deg_b(\Gamma_1\circ\Gamma_2)=deg_b(\Gamma_1)+deg_b(\Gamma_2).$$
The above composition does not invary the class of Lie admissible graphs
corresponding to linear Poisson structures.
Since we are interested in the graphs not in the kernel of the
Kontsevich representation,
we will consider the truncation of the above composition
due to the (orthogonal) projection $Pr$ from all admissible graphs
to our class of Lie admissible edge-labeled graphs $\tilde{\C{G}}$.
The resulting composition of graphs is now an internal operation, 
still graded by $deg_b$.
\begin{defin}\label{D:star1}
The {\em internal composition of graphs} of $\tilde{\C{G}}$ is defined as follows:
$$\tg_1\circ\tg_2=Pr[\sum_{i=1}^m (-1)^{(i-1)(m'-1)}\tg_1\circ_i\tg_2],
\quad \tg_1\in \tilde{\C{G}}_{n,m}, \ \tg_2\in \tilde{\C{G}}_{n',m'},$$
where $\circ_i$ is the insertion of $\tg_2$ at the $i^{th}$
boundary vertex of $\tg_1$ using ``Leibniz rule''
i.e. summing over all possible graphs where the ``$i^{th}$ legs'' of $\tg_1$ lend on 
vertices of $\Gamma_2$, internal and external.
The edge-labeling of the resulting graph is inherited from 
the edge-labeling of the two graphs $\tg_i$.
\end{defin}
Graph composition is compatible with the equivalence relations $ng^2$ and $ng$ (\cite{Po}, p.5).
\begin{lem}
There is a unique graph composition on $\C{G}$ and on $H$ such that
the canonical projections $ng^2$ and $ng$ are morphisms of pre-Lie algebras.
As a consequence, taking the equivalence class of a graph $[\ ]$ 
is also a pre-Lie algebra morphism.
\end{lem}
If $I_i$ denotes the set of incoming edges at the $i^{th}$ boundary point
of $\Gamma_1$ and $[n_2],[m_2]$ denote the sets of 
internal and respectively external vertices of $\Gamma_2$,
then the ``Leibniz rule'' at the $i^{th}$ vertex yields:
$$\Gamma_1\circ_i\Gamma_2=\sum_{f:I_i\to [n_2]\cup[m_2]}\Gamma_1\circ_i^f\Gamma_2,$$
where $\circ_i^f$ denotes the operation of 
replacing the vertex $i$ with a disjoint copy of the set $[n']$, 
with the edges $e=(v\to i)\in I_i$ now pointing to $f(e)$.
The ``two components'' of $f$, $f_i, f_b$ denote their co-restriction
to internal and boundary vertices, respectively.

Now in order for the resulting graph to have internal vertices with only one incoming arrow,
the component $f_i$ must be injective,
yielding the following formula for graph composition.
\begin{lem}\label{L:comb1}
If $\Gamma_1\in\C{G}_{n_1,m_1}$ and $\Gamma_2\in \C{G}_{n_2,m_2}$ then 
$$\Gamma_1\circ_i\Gamma_2=
\sum_{f_i\cup f_b:I_i\hookrightarrow [n_2]\cup[m_2]}\Gamma_1\circ_i^f\Gamma_2,$$
where $f_i$ is 1:1.
\end{lem}

We are now ready to prove that the ``sum of all graphs'' is a solution.


\def\Vin{V^{\rm in}}
\def\Vbd{V^{\rm bd}}

\section{The canonical solution}\label{S:solution}
Instead of investigating whether the Kontsevich solution 
restricts to graphs without circuits to a cocycle of the corresponding dg-coalgebra \cite{FI},
still providing a star-product,
we will provide a canonical solution of the deformation equation (see also \cite{AIS}).

Let $Z_k=\sum_{\Theta\in G_{k,2}}\Theta/|Aut(\Theta)|$ and $Z=\sum Z_k h^k$.
In order to prove $Z\circ Z=0$ we need to investigate the coefficients of
$$\sum_{i+j=n;\ i,j\geq 0}Z_i\circ_1 Z_j-Z_i\circ_2 Z_j=
\sum_{\Gamma\in G_{n,3}} B_{\Gamma}\cdot\Gamma,$$
where the coefficient $B_\Gamma$ is the difference between 
the coefficients (possibly zero) of the graph $\Gamma$ 
resulting from left and from right graph insertions ($\circ_1$ and $\circ_2$):
$$B_\Gamma=B^L_\Gamma-B^R_\Gamma.$$
To simplify notation, for any $\Gamma\in G_{n,m}$,
$\underline{\Gamma}$ denotes the corresponding normalized basis element,
i.e. $\Gamma/|Aut(\Gamma)|$.
The {\em normalized bases} of $kG$ is: $\{\underline{\Gamma}\}_{\Gamma\in G_{n,m}}$.

The key fact (Proposition \ref{P:inj}) is that $\circ_1$ is ``injective'' 
(similarly $\circ_2$), i.e. from the composition $\Gamma_1\circ_1\Gamma_2$ 
one can recover the operands $\Gamma_1$ and
$\Gamma_2$ (``left groupoid structure'' $\Gamma:\Gamma_1\to \Gamma_2$).
In general the pair $(\Gamma_1,\Gamma_2)$ responsible for a summand $\Gamma$
as a result of a left insertion $\circ_1$ is different from the unique pair
yielding a sum involving $\Gamma$ in a right insertion $\circ_2$
(Is the ``left groupoid'' isomorphic to the ``right groupoid''?).

Now comparing the two sums:
\begin{equation}\label{E:sigma}
\Sigma_k=\sum_{i+j=n, i,j\ge0}\quad \sum_{\Gamma_1\in\C{G}_{i,2},\Gamma_2\in \C{G}_{j,2}}
\Gamma_1\circ_k\Gamma_2,\quad k=1,2
\end{equation}
corresponding to left, and respectively right insertions,
one obtains that the respective coefficients are equal (Corollary \ref{cor:id}), 
a fact expected due to the left/right symmetry,
and proved  as The Coefficient Theorem \ref{conj}.
\begin{th}\label{T:main}
If $Z=\sum_{\Theta\in G_{n,2}}\underline{\Theta}$ then  $[Z,Z]=0$.
\end{th}
To prove the above claims, we start with some preparatory lemmas.

For $\Gamma\in\C{G}$, let $\Vin(\Gamma)$ denote its set of internal 
vertices, $\Vbd(\Gamma)$ its set of boundary vertices, and let $V(\Gamma)=\Vin\cup \Vbd$.
For $u\in \Vin(\Gamma)$, let $u_L$ and $u_R$ be the left and right 
descendants of $u$, respectively. 
Moreover, denote by $(u,\ldots,v)$ 
a directed path starting at $u$ and ending at the vertex $v$.
\begin{lem}\label{lem:LRpath}
Let $\Theta\in G_{n,2}$ with boundary vertices $L$ and $R$. 
Then for each internal vertex $u$ of $\Theta$, there 
is a directed path from $u$ to $L$ and a directed path from $u$ to $R$
\end{lem}
\begin{pf}
Define a partial order on internal vertices corresponding to the ``flow'' direction
corresponding to the oriented edges (no loops!).
Since any internal vertex has two descendants, clearly there is a path starting at $u$ ending
at a boundary point, say $L$.
Now not all paths may end at $L$, since one may trace back last arrow and descend on the other arrow,
until the end of the path is not $L$. 

In particular, binary graphs without loops are connected.
\end{pf}
\begin{rem}
Note that the lemma may fail for graphs with loops,
and for $m=0,1$, $\C{G}_{n,m}$ contains no binary admissible graph without loops.
\end{rem}
\begin{lem}\label{lem:norm}
Let $\Gamma\in G_{n,3}$ and $\Theta$ be a normal subgraph of $\Gamma$,
i.e. $\Gamma/\Theta$ is still admissible,
with at least one boundary point.

\n $(i)$ If $u\in \Vin(\Theta)$ then $u_L,u_R\in V(\Theta)$
$(\Theta$ is a ``total subgraph'' of $\Gamma)$.
 
\n $(ii)$ 
If $(v_1,\ldots, v_t)$ is a path from an interior point $v_1$ of $\Theta$
to a boundary point $v_t$ of $\Gamma$ then 
$v_i\in V(\Theta)$ for all $1\leq i\leq t$
$(\Theta$ is ``geodesically complete''$)$.
\end{lem}
\begin{pf}
$(i)$ follows from the fact that if $u\in V(\Theta)$ and say
$u_L\not\in V(\Theta)$,
then the edge $(u,u_L)$ is not present in $\Theta$. 
This would contradict that $\Theta$ is a binary graph,
since it has an internal vertex $u$ with at most 
one outgoing edge. 

$(ii)$ follows from a recursively application of $(i)$,
using Lemma \ref{lem:LRpath}.
\end{pf}
Our next goal is to prove that for each graph $\Gamma\in\C{G}_{n,3}$
there is a unique factorization in terms of graphs with less boundary points:
$$\Gamma=\Gamma_1\circ_1\Gamma_2\quad (\Gamma=\Gamma_1'\circ_2\Gamma_2').$$
Each such decomposition will correspond to a ``maximal factor'' of $\Gamma$,
so here too ``maximal implies prime''!
\begin{lem}
For $\Gamma\in G_{n,3}$ there are unique normal subgraphs of $\Gamma$,
denoted $\alpha_L(\Gamma), \alpha_R(\Gamma)\in G_{n,2}$,
sitting on the leftmost, and respectively rightmost, 
two boundary vertices of $\Gamma$.
\end{lem}
\begin{pf}
Recall that being normal ensures that the quotient
$\Gamma/\alpha_L(\Gamma)$ is still a binary (exactly two descendants) admissible graph.

Suppose that $\Gamma_1$ and $\Gamma_2$ are two different normal subgraphs of 
$\Gamma$ sitting on say boundary points $1$ and $2$ of $\Gamma$
(the other case follows by symmetry)
$$b^L_0=\{1,2\}\subset\Gamma_i\subset \Gamma.$$
Then there exist an internal vertex $u$ of $\Gamma_1$ but not in $\Gamma_2$,
since they cannot both equal $b^L_0$.
Note that by Lemma \ref{lem:norm}, 
any path starting at $u$ must end at a boundary vertex: $1$ or $2$.

Since $\Gamma_2$ is normal, $\Gamma_2'=\Gamma/\Gamma_2\in G_{n,2}$ 
is a binary admissible graph. 
By definition, we have $u\in \Vin(\Gamma_2')$. 
However, there is no directed path from $u$ to the right boundary vertex $\Gamma_2'$,
contradicting Lemma~\ref{lem:LRpath}.
\end{pf}
\begin{defin}
For $\Gamma\in G_{n,3}$, 

\bs\n $(i)$ $E_\Gamma^L=\underline{(\Gamma/\alpha_L(\Gamma))}\circ_1\underline{\alpha_L(\Gamma)}$
and $E_\Gamma^R=\underline{(\Gamma/\alpha_R(\Gamma))}\circ_2\underline{\alpha_R(\Gamma)}$.

\bs\n $(ii)$ $C_\Gamma^L=<E_\Gamma^L,\Gamma>$ and 
$C_\Gamma^R=<E_\Gamma^R,\Gamma>$ are the coefficients of $\Gamma$ in $E_\Gamma^L$, 
and $E_\Gamma^R$ respectively.

\bs\n $(iii)$ $C_\Gamma=C_\Gamma^L-C_\Gamma^R$.
\end{defin}
We now prove the key fact, that the left and right insertions are ``injective''.
\begin{prop}\label{P:inj}
Let $\Gamma',\Gamma''\in \C{G}_{\bullet,2}$ and $\Gamma\in\C{G}_{n,3}$.

(i) If $<\Gamma''\circ_1\Gamma',\Gamma>\ne 0$ then
$\Gamma'=\alpha_L(\Gamma), \quad \Gamma''=\Gamma/\alpha_L(\Gamma)$.

(ii) If $<\Gamma''\circ_2\Gamma',\Gamma>\ne 0$ then
$\Gamma'=\alpha_R(\Gamma), \quad \Gamma''=\Gamma/\alpha_R(\Gamma)$.
\end{prop}
\begin{pf}
Let $\Gamma',\Gamma''\in G_{n,2}$ be such that $\Gamma\in G_{n,3}$ is a summand 
of $\Gamma''\circ_1\Gamma'$. 
Then, there exists a way to land legs from the 
left boundary vertex of $\Gamma''$ onto the vertices of $\Gamma'$ so that the 
resulting graph is $\Gamma$.

Since $\Gamma'$ is a normal subgraph of $\Gamma$ sitting on the 
its 1st and 2nd boundary vertex, it follows that $\Gamma'\subseteq \alpha_L(\Gamma)$.
Moreover $\Gamma'$ is the maximal subgraph of $\Gamma$ sitting on its 1st and 2nd 
boundary vertex. 
For otherwise, $\alpha_L(\Gamma)$ contains an internal vertex $u$ of $\Gamma''$. 
Now, by Lemma~\ref{lem:LRpath}, there exist a path $P_u^R$ from $u$ to the second 
boundary vertex of $\Gamma''$. 
By Lemma~\ref{lem:norm}, all 
the vertices (internal and external) in $P_u^R$ are in $\alpha_L(\Gamma)$. In particular, 
the second boundary vertex of $\Gamma''$ (which is the third boundary vertex of $\Gamma$) 
would also have to be in  $\alpha_L(\Gamma)$. This would contradict the fact that $\alpha_L(\Gamma)$ 
has to sit on the  1st and 2nd boundary vertex of $\Gamma$. Thus $\Gamma'=\alpha_L(\Gamma).$

Since $\Gamma$ is a summand of $\Gamma''\circ_1\Gamma'$ and $\Gamma'=\alpha_L(\Gamma)$, we 
have $\Gamma''=\Gamma/\alpha_L(\Gamma)$; because $\circ_1$ splits the edges landing on the 
1st boundary vertex of $\Gamma''$ and land them on vertices of $\Gamma'=\alpha_L(\Gamma)$ 
while collapsing $\alpha_L(\Gamma)$ in $\Gamma$ does the converse,
recovering $\Gamma''$.

Similarly, if $\Gamma',\Gamma''\in G_{n,2}$ are such that $\Gamma$ is a summand of 
$\Gamma''\circ_2\Gamma'$ then $\Gamma''=\alpha_R(\Gamma)$ and 
$\Gamma'=\Gamma/\alpha_R(\Gamma)$.
\end{pf}
Regarding graph insertions as partially defined binary operations,
the above result may be rephrased as follows.
\begin{cor}\label{C:UF}
Boundary graph insertions have the unique factorization property.
\end{cor}
As an immediate consequence we obtain that the $\Gamma$-coefficients of $[Z,Z]$
result from a unique left/right composition,
namely the composition of the unique normal maximal left/right suported subgraphs.
\begin{cor}\label{lem:B=C}
$B_\Gamma=C_{\Gamma}$.
\end{cor}
\begin{pf}
As a consequence of Proposition~\ref{P:inj}, 
$B_\Gamma=B^L_\Gamma-B^R_\Gamma$
represents the contributions from a left composition
of a unique pair of graphs $(\Gamma_1,\Gamma_2)$
and of a right composition of a unique pair $(\Gamma_1',\Gamma_2')$. 
The corresponding multiplicities are $C^L_\Gamma$ and $C^R_\Gamma$.
Therefore $B_\Gamma=C_\Gamma$.
\end{pf}
All that is left in order to prove the main theorem,
is to prove that left insertions produce the same coefficients as right insertions,
i.e. $C^L_\Gamma=C^R_\Gamma$.
Fix a summand $\Gamma$ of a fixed pair of graphs $\Gamma_1, \Gamma_2$,
i.e. $\Gamma$ has a non-trivial coefficient $C^L_\Gamma$ in the sum expressing 
the left boundary composition $\Gamma_1\circ_1\Gamma_2$.
Then there is a {\em left extension}:
$$\Gamma_2\overset{\pi}{\to}\Gamma\to\Gamma_1$$
characterized by the {\em insertion data} $\pi$ 
(see section \ref{S:proof} for additional details),
with $\Gamma_2$ collapsing to the left boundary vertex of $\Gamma_1$.
%
\begin{conj}\label{conj}
After normalization, the non-trivial coefficient of $\Gamma$ 
as a summand of the left insertion operation of $\Gamma_2$ in $\Gamma_1$
is 1.
Therefore, if non-trivial, the left/right normalized multiplicities are:
$$\underline{L}^\Gamma_{\Gamma_1\Gamma_2}=
<\underline{\Gamma_1}\circ_1\underline{\Gamma_2},\underline{\Gamma}>
=1=<\underline{\Gamma_1}\circ_1\underline{\Gamma_2},\underline{\Gamma}>
=\underline{R}^\Gamma_{\Gamma_1\Gamma_2},$$
where $\pi_L,\pi_R$ are the left/right insertion data
determined by the graph $\Gamma$.
\end{conj}
We will first exploit the result, 
defering the proof to section \ref{S:proof}.
\begin{cor}\label{cor:id}\

(i) For all $\Gamma\in\C{G}_{n,3}$,
its left multiplicity equals its right multiplicity:
$$C^L_\Gamma=C^R_\Gamma.$$

(ii) $\forall \Gamma\in\C{G}_{\bullet,3}$:
$$ <\Sigma_1,\Gamma>=<\Sigma_2,\Gamma>.$$

\end{cor}
\begin{pf}
Any $\Gamma\in \C{G}_{n,3}$ appears as part of $(\Gamma/\alpha_L(\Gamma))\circ_1\alpha(\Gamma)$.
The coefficient of $\Gamma$ in both $\Sigma_1$ and $\Sigma_2$ (Equation \ref{E:sigma}) is 
$|Aut(\Gamma)|$, i.e (i) holds, and the two sums are equal.
\end{pf}
>From Lemma~\ref{lem:B=C} $B_\Gamma=C_\Gamma$,
which by the Corollary~\ref{cor:id} vanish for all graphs $\Gamma$.
This implies that 
$$\sum_{i+j=n;\ i,j\geq 0}Z_i\circ Z_j=0,$$
which yields $Z\circ Z=0$,
concluding the proof of the Main Theorem \ref{T:main}.

%
%
\subsection{Examples}
Consider the graphs $\Gamma_1,\Gamma_2,\Gamma_3\in \C{G}_{1,3}$, defined as follows:
$$
\begin{xy}
,(-10,4)*{\Gamma_1=\quad };
,(5,7)*{\bullet};(5,-4)*{\circ}**\dir{-} ?>*\dir{>}
,(5,7);(15,-4)*{\circ}**\dir{-}?>*\dir{>}
,(15,7);(5,-4)*{\circ}
,(15,-4)*{\circ}
,(-5,-4)*{\circ};
\end{xy}
\ \   
\begin{xy}
,(-12,4)*{\mbox{ , }\Gamma_2=\quad };
,(5,7)*{\bullet};(-5,-4)*{\circ}**\dir{-} ?>*\dir{>}
,(5,7);(15,-4)*{\circ}**\dir{-}?>*\dir{>}
,(5,-4)*{\circ};
\end{xy}
\ \   
\begin{xy}
,(-6,4)*{,\mbox{ and }\Gamma_3=\ }; (5,-4)*{\circ}
,(5,7);(15,-4)*{\circ}
,(15,7);(5,-4)*{\circ}**\dir{-} ?>*\dir{>}
,(15,7)*{\bullet};(15,-4)*{\circ}**\dir{-} ?>*\dir{>}
,(25,-4)*{\circ};
\end{xy}
$$

(1) {\bf The constant Case.} 
Any admissible graphs $\Gamma\in G_{n,3}$ can be expressed as 
 $\Gamma=\Gamma_1^r\Gamma_2^s\Gamma_3^t$, where $\Gamma_i$, $i=1,2,3$, is 
as defined earlier.
Thus  
$$E_\Gamma=\underline{{b_1^{r+s}}}\circ_1\underline{{b_1^t}}
-\underline{{b_1^{s+t}}}\circ_2\underline{{b_1^r}}.$$
Since the coefficient of $\Gamma$ in $\underline{{b_1^{r+s}}}\circ_1\underline{{b_1^t}}$
is $\binom{r+s}{s}$ and  the coefficient of $\Gamma$ in 
$\underline{{b_1^{s+t}}}\circ_2\underline{{b_1^r}}$ is $\binom{s+t}{s}$ then 
$$C_{\Gamma}=\frac{\binom{r+s}{s}}{t!(r+s)!}-\frac{\binom{s+t}{s}}{r!(s+t)!}=0$$
In the normalized bases $\underline{C}^L_{\Gamma}=\underline{C}^R_{\Gamma}=1$
(similarly for the other normalized coefficients bellow).

(2) {\bf The linear case with $n=2$.}
There are $9$ admissible graphs in $G_{2,3}$: $\Gamma_1^2$, $\Gamma_2^2$, $\Gamma_3^2$, 
$\Gamma_1\Gamma_3$, $t_2^L$, $t_2^R$ $c_2^L$, $c_2^R$, and $c_2$.
\begin{enumerate}
\item Let $\Gamma=\Gamma_1^2$. We have $\alpha_L(\Gamma)=b_0$, $\alpha_R(\Gamma)=b_1^2$, 
$\Gamma/\alpha_L(\Gamma)=b_1^2$, and $\Gamma/\alpha_R(\Gamma)=b_0$. Thus,
the coefficient of $\Gamma$ in 
$\underline{b_1^2}\circ_1 b_0-\underline{b_0}\circ_2 \underline{b_1^2}$,
which is $C_\Gamma=\frac12-\frac12=0$. 
\item Let $\Gamma=\Gamma_2^2$. We have $\alpha_L(\Gamma)=b_1^2$, $\alpha_R(\Gamma)=b_0$, 
$\Gamma/\alpha_L(\Gamma)=b_0$, and $\Gamma/\alpha_R(\Gamma)=b_1^2$. Thus,
the coefficient of $\Gamma$ in 
$b_0\circ_1 \underline{b_1^2}-\underline{b_1^2}\circ_2 \underline{b_0}$,
which is $C_\Gamma=\frac12-\frac12=0$.
\item Let $\Gamma=\Gamma_3^2$. We have $\alpha_L(\Gamma)=b_0$, $\alpha_R(\Gamma)=b_0$, 
$\Gamma/\alpha_L(\Gamma)=b_1^2$, and $\Gamma/\alpha_R(\Gamma)=b_1^2$. Thus,
the coefficient of $\Gamma$ in 
$\underline{b_1^2}\circ_1 b_0-\underline{b_1^2}\circ_2 \underline{b_0}$,
which is $C_\Gamma=\frac12-\frac12=0$.
\item Let $\Gamma=\Gamma_1\Gamma_3$. We have $\alpha_L(\Gamma)=b_1$, $\alpha_R(\Gamma)=b_1$, 
$\Gamma/\alpha_L(\Gamma)=b_1$, and $\Gamma/\alpha_R(\Gamma)=b_1$. Thus,
the coefficient of $\Gamma$ in   
$\underline{b_1}\circ_1 b_1-\underline{b_1}\circ_2 \underline{b_1}$,
which is $C_\Gamma=1-1=0$.
\item Let $\Gamma=t_2^L$. We have $\alpha_L(\Gamma)=b_0$, $\alpha_R(\Gamma)=b_1$, 
$\Gamma/\alpha_L(\Gamma)=b_2^L$, and $\Gamma/\alpha_R(\Gamma)=b_1$. Thus,
the coefficient of $\Gamma$ in 
$\underline{b_2^L}\circ_1 \underline{b_0}-\underline{b_1}\circ_2 \underline{b_1}$,
which is $C_\Gamma=1-1=0$.
\item Let $\Gamma=t_2^R$. We have $\alpha_L(\Gamma)=b_1$, $\alpha_R(\Gamma)=b_0$, 
$\Gamma/\alpha_L(\Gamma)=b_1$, and $\Gamma/\alpha_R(\Gamma)=b_2^R$. Thus,
the coefficient of $\Gamma$ in 
$\underline{b_1}\circ_1 \underline{b_1}-\underline{b_2^R}\circ_2\underline{b_0}$,
which is $C_\Gamma=1-1=0$.
\item Let $\Gamma=c_2^L$. We have $\alpha_L(\Gamma)=b_1$, $\alpha_R(\Gamma)=b_0$, 
$\Gamma/\alpha_L(\Gamma)=b_1$, and $\Gamma/\alpha_R(\Gamma)=b_1^2$. Thus,
the coefficient of $\Gamma$ in 
$\underline{b_1}\circ_1 \underline{b_1}-\underline{b_1^2}\circ_2\underline{b_0}$,
which is $C_\Gamma=1-\frac{2}{2}=0$.
\item Let $\Gamma=c_2^R$. We have $\alpha_L(\Gamma)=b_0$, $\alpha_R(\Gamma)=b_1$, 
$\Gamma/\alpha_L(\Gamma)=b_1^2$, and $\Gamma/\alpha_R(\Gamma)=b_1$. Thus,
the coefficient of $\Gamma$ in 
$\underline{b_1^2}\circ_1\underline{b_0}-\underline{b_1}\circ_2 \underline{b_1}$,
which is $C_\Gamma=\frac{2}{2}-1=0$.
\item Let $\Gamma=c_2$. We have $\alpha_L(\Gamma)=b_0$, $\alpha_R(\Gamma)=b_0$, 
$\Gamma/\alpha_L(\Gamma)=b_2^R$, and $\Gamma/\alpha_R(\Gamma)=b_2^L$. Thus,
 the coefficient of $\Gamma$ in 
$\underline{b_2^R}\circ_1\underline{b_0}-\underline{b_2^L}\circ_2 \underline{b_0}$,
which is $C_\Gamma=1-1=0$.
\end{enumerate}
Hence $\sum_{i+j=2; i,j\geq 0}[Z_i,Z_j]=\sum_{\Gamma\in G_{2,3}} C_{\Gamma}\cdot\Gamma=0.$

%
%
\subsection{Proof of Theorem \ref{conj}}\label{S:proof}
We prove that the multiplicity of a graph as a summand in a 
graph composition is only due to their groups of symmetries.

Fix graphs $\Gamma_1, \Gamma_2$ and a summand $\Gamma$ of their left boundary insertion:
$$L_{\Gamma_1\Gamma_2}^{\Gamma}=<\Gamma_1\circ_1 \Gamma_2,\Gamma>\ne 0.$$
Since similar considerations apply to right insertions
and to the corresponding coefficient $R^\Gamma_{\Gamma_1\Gamma_2}$,
we will use the generic notation $C^\Gamma_{\Gamma_1\Gamma_2}$.

Then there is a left graph extension:
$$\Gamma_2\overset{\pi}{\to}\Gamma\to \Gamma_1,$$
determined by the left insertion data $\pi:S\subset V_1 \to T\subset V_2$
defining the way the left leg arrows of $\Gamma_1$
land on the vertices of $\Gamma_2$, internal or boundary.
Each insertion data $\pi$ yields an admissible graph $\Gamma_\pi$. 
Its isomorphism class will be called the {\em type of the insertion}.

Recall that for a linear Poisson structure, 
the non-boundary portion of the insertion $\pi$ is injective.

Let $\C{D}$ be the set of all insertion data $\pi$.
For any $\pi\in\C{D}$, 
let $\C{D}_{\Gamma}\subseteq \C{D}$ be those insertion data 
of the same type as $\Gamma$. 
Then 
\begin{eqnarray}\label{CfD}
C_{\Gamma_1\Gamma_2}^{\Gamma}=<\Gamma_1\circ_1 \Gamma_2,\Gamma>=|\C{D}_{\Gamma}|.
\end{eqnarray}
For any insertion data $\pi$, let $Aut(\Gamma_\pi,\pi)$ be the set of automorphism 
$Aut(\Gamma_\pi)$ that fix $\pi$.

We claim that the multiplicity of a summand in a left (right) boundary composition
is given by the following formula.
\begin{lem}\label{C_Gamma=}
$$|\C{D}_{\Gamma}|=\frac{|Aut(\Gamma_1)|\cdot|Aut(\Gamma_2)|}{|Aut(\Gamma_\pi,\pi)|}$$
\end{lem}
We delay the proof of Lemma~\ref{C_Gamma=} to make some general observations.

Consider the action $\tau$ of $H=Aut(\Gamma_1)\times Aut(\Gamma_2)$ on $\C{D}_\Gamma$ 
defined as follows.
For all $\rho=(\rho_1,\rho_2)\in H$ and for all $\pi\in \C{D}$
with $\pi:S\subset V_1 \to T_\pi\subset V_2$, we have 
\begin{eqnarray}\label{def_action}
\tau(\rho,\pi)=\pi_\rho:\rho_1(S)=S&\to& \rho_2(T_\pi)\\
                                x&\mapsto&\rho_2\pi\rho_1^{-1}(x).
\end{eqnarray}
\begin{cl}\label{orb}
For $\pi\in\C{D}$, we have
$$\C{O}(\pi)=\left\{\pi_\rho:\ \rho\in H \right\}=\C{D}_{\Gamma_\pi},$$ 
i.e. the action action $\tau$ is transitive on $\C{D}_{\Gamma_\pi}$.
\end{cl}
\begin{pf}
By definition, $\C{O}(\pi)\subseteq \C{D}_{\Gamma_\pi}$. Conversely, 
if $\Gamma_{\pi'}\in \C{D}_{\Gamma_\pi}$ then there exits $\pi'\in\C{D}$ 
such that $\Gamma_{\pi'}\cong\Gamma_\pi$. Hence, there exist 
$\phi\in Aut(\Gamma_\pi)$ such that $\phi(\C{D}_{\Gamma_\pi})=\Gamma_{\pi'}$.
Now let $\rho_\phi=(\phi_{|\Gamma_1},\phi_{|\Gamma_2})\in H$; then 
$\rho_\phi(\pi)=\pi'\in \C{O}(\pi)$. Thus $\C{D}_{\Gamma_\pi}\subseteq\C{O}(\pi)$,
and the claim follows.
\end{pf}
\begin{cl}\label{stab}
For $\pi\in\C{D}$, we have 
$$|Stab(\pi)|=|\left\{\rho\in H:\ \pi_\rho=\pi\right\}|=|Aut(\Gamma_\pi,\pi)|.$$
\end{cl}
\begin{pf}
We show that there exist a bijection
$f: Aut(\Gamma_\pi,\pi)\to Stab(\pi)$. 
For $\phi\in Aut(\Gamma_\pi)$, define 
$$f(\phi)=(\rho_1,\rho_2),$$
by first restricting $\phi$ to the unique normal subgraph $\Gamma_1$, which 
therefore is invaried by 
$$\rho_1(\phi)=\phi_{|\Gamma_1}.$$
Then, $\phi$ induces an automorphism of the quotient:
$$\rho_2(\phi)=\phi_{|(\Gamma/\Gamma_1)}=\phi_{|\Gamma_2}.$$
Thus, $f(\phi)=(\rho_1,\rho_2)\in Stab(\pi)$ since, by definition of $\rho_1$ and $\rho_2$., 
$\rho_2\pi\rho_1^{-1}=\pi$. 
. 
It is easy to see that $f$ is injective, since $V=V_1\cup V_2$.

To prove that $f$ is surjective, let $\rho=(\rho_1,\rho_2)\in Stab(\pi)$. 
Then there exist unique automorphisms $\phi_1,\phi_2\in Aut(\Gamma_\pi)$ obtained 
by extending $\rho_1$ and $\rho_2$ in such a way that 
${\phi_1}_{|\Gamma_2}=id_{\Gamma_2}$ and ${\phi_2}_{|\Gamma_1}=id_{\Gamma_1}$. 
Thus $\phi=\phi_1\phi_2\in Aut(\Gamma_\pi,\pi)$ is such that 
$\phi_{|\Gamma_1}=\rho_1$ and $\phi_{|\Gamma_2}=\rho_2$, i.e. $f(\phi)=(\rho_1,\rho_2)$.
Since $f$ is a bijection, we have $|Stab(\pi)|=|Aut(\Gamma_\pi,\pi)|$.
\end{pf}

\begin{pf}[Proof of Lemma~\ref{C_Gamma=}]\

Using an orbit-stabilizer argument, \eqref{CfD}, and Claims~\ref{orb}\&\ref{stab}, 
we obtain
\begin{eqnarray}\label{orb_size}
C_{\Gamma_1\Gamma_2}^{\Gamma}=
\C{D}_{\Gamma_\pi}=|\C{O}(\pi)|=\frac{|H|}{|Stab(\pi)|}=
\frac{|Aut(\Gamma_1)|\cdot|Aut(\Gamma_2)|}{|Aut(\Gamma_\pi,\pi)|},
\end{eqnarray}
proving Lemma~\ref{C_Gamma=}. 
\end{pf}
\begin{pf}[Proof of The Coefficient Theorem~\ref{conj}]\

For $\Gamma\in G_{n,3}$, let $\Gamma_1^L=\alpha_L(\Gamma)$,  $\Gamma_2^L=
\Gamma/\Gamma_1^L$, $\Gamma_1^R=\alpha_R(\Gamma)$, and  $\Gamma_2^R=
\Gamma/\Gamma_1^R$. Recall that 

$C_\Gamma^L=<\underline{\Gamma}_1^L\circ_1\underline{\Gamma}_2^L,\Gamma>$,
$C_\Gamma^R=<\underline{\Gamma}_1^R\circ_2\underline{\Gamma}_2^R,\Gamma>$, and  
$C_\Gamma=C_\Gamma^L-C_\Gamma^R.$

\bs If $C_\Gamma\ne 0$ then $<\Gamma_1^L\circ_1\Gamma_2^L,\Gamma>\ne 0$ and 
$<\Gamma_1^R\circ_2\Gamma_2^R,\Gamma>\ne0$. 
Hence,  there exist two insertion data $\pi_L$ and $\pi_R$ such that 
$\Gamma_{\pi_L}\cong\Gamma\cong\Gamma_{\pi_R}$.
Moreover, it follows from Lemma~\ref{C_Gamma=} that 
$$L_{\Gamma_1^L,\Gamma_2^L}^{\Gamma}=|\C{D}_{\Gamma}|=
\frac{|Aut(\Gamma_1^L)|\cdot|Aut(\Gamma_2^L)|}{|Aut(\Gamma,\pi_L)|},$$
$$R_{\Gamma_1^R,\Gamma_2^R}^{\Gamma}=|\C{D}_{\Gamma}|=
\frac{|Aut(\Gamma_1^R)|\cdot|Aut(\Gamma_2^R)|}{|Aut(\Gamma,\pi_R)|}.$$
Thus 
$$C_\Gamma^L=\frac{L_{\Gamma_1^L,\Gamma_2^L}^{\Gamma}}
{|Aut(\Gamma_1^L)|\cdot|Aut(\Gamma_2^L)|}=
\frac{1}{|Aut(\Gamma,\pi_L)|},$$
$$C_\Gamma^R=\frac{R_{\Gamma_1^R,\Gamma_2^R}^{\Gamma}}
{|Aut(\Gamma_1^R)|\cdot|Aut(\Gamma_2^R)|}=
\frac{1}{|Aut(\Gamma,\pi_R)|}.$$
Now it remains to show that if $\Gamma_{\pi_L}\cong\Gamma_{\pi_R}$ then
$$|Aut(\Gamma,\pi_L)|=|Aut(\Gamma,\pi_R)|.$$
In fact both automorphism groups equal $Aut(\Gamma)$.
In order to prove this, 
note that there are (natural) restriction monomorphisms:
$$Aut(\Gamma^L_1)\times Aut(\Gamma^L_2)\leftarrow Aut(\Gamma)
\to Aut(\Gamma^R_1)\times Aut(\Gamma^R_2),$$
since $\Gamma^L_i$, and respectively $\Gamma^R_i$,
are invaried as being (unique) maximal left normal subgraphs.

\begin{lem}
$$Aut(\Gamma,S)=Aut(\Gamma,\pi),$$
where $S$ is the domain of $\pi$ and $Aut(\Gamma,S)$
is the subset of automorphisms of $\Gamma$ which invary $S$,
i.e. $\Phi(S)\subset S$.
\end{lem}
\begin{pf}
It is enough to prove $``\subset''$, since the other inclusion
follows from the definition of $Aut(\Gamma,\pi)$.

If $\Phi\in Aut(\Gamma)$ and $\Phi(S)\subset S$ then 
$$\Phi(s\to t)=\Phi(s)\to \Phi(t)=s'\to \Phi(t).$$
Since $S$ has the property that any of its points 
has a unique arrow towards $V_2$, the vertices of $\Gamma_2$,
then $\Phi(t)=\pi(s')=\pi(\Phi(s))$,
i.e. $\pi\Phi=\Phi\pi$ on $S$,
and therefore $\Phi$ invaries $\pi$.
\end{pf}
Now the unique factorization implies that the ``Galois group''
$Aut(\Gamma,\pi)$ is the full automorphism group.
\begin{lem}
$$Aut(\Gamma,\pi)=Aut(\Gamma).$$
\end{lem}
\begin{pf}
Let $\pi$ be the left insertion data yielding $\Gamma$ as a left extension
(by unique factorization).
If $\Phi\in Aut(\Gamma)$ not only $\Phi$ invaries $\Gamma_1$ and $\Gamma_2$,
but also $S$, the domain of $\pi$ as being the set of arrows lending on 
the left leg of $\Gamma_1$.

By the previous lemma, $\Phi$ invaries $\pi$.
\end{pf}
Therefore $Aut(\Gamma,\pi_L)=Aut(\Gamma)=Aut(\Gamma,\pi_R)$
is the satbilizer of the action and the normalized coefficients are trivial
or equal to 1.

This concluds the proof of Theorem \ref{conj}.
\end{pf}
\begin{rem}\label{R:linear} 
In the case of a linear Poisson structures the structure of the Galois group
of a left (right) extension is simpler 
(subobject of the fibered product of $Aut(\Gamma_1)$ and $Aut(\Gamma_2)$),
since $\pi$, the insertion data, is injective at the level of interior points,
and a permutation of $S$ is equivalent to an inverse permutation of $T$.

Nevertheless the ``simplification'' entailing the left-right symmetry (equal coefficients)
seams to be due to the lack of circuits, rather than, as one might expect,
from the one-incoming arrow property satisfied by graphs in the linear case 
($\pi$ injective on interior points).
\end{rem}

\section{Conclusions and further developments}\label{S:conclusions}
We proved the existence of solutions of Maurer-Cartan equation, 
without assuming the Jacobi identity holds,
implying that the primary obstruction to have a full deformation vanishes anyway
(the cohomology class).
The main fact used in the proof is the unique factorization property enjoied by graph insertions.

As a first problem to be investigated we note the question regarding the relation 
with the other ``universal solution'', the Haussdorf series, 
living on the base space.
It is also natural to look for a physical interpretation of our solution
as a (semi-classical part of the) correlation function in the spirit of \cite{CF}.

We believe that these are interesting topics for further study,
revealing some of the intimate relashionship between 
the mathematics and physics of quantum phenomena.


\end{document}